\documentclass[a4paper]{amsart}
\usepackage{epsfig}
\usepackage{amsmath}
\usepackage{amssymb}
\usepackage{bbm}
\usepackage{a4wide}
\parskip\smallskipamount

\let\set\mathbbm
\def\<#1>{\langle#1\rangle}

\newcommand\todo[1][.]{\edef\tmpa{.}\edef\tmpb{#1}%
  \ifx\tmpa\tmpb
    \typeout{To Be on page \thepage}\fbox{\bf To Be}
  \else
    \typeout{To Be on page \thepage: #1}\fbox{{\bf To Be:} #1}
  \fi
}

\newtheorem{ethm}{Empirical Result}

\begin{document}

 \title{Experiments with a Positivity Preserving Operator}

 \author[Manuel Kauers]{Manuel Kauers\,$^\ast$}
 \address{Manuel Kauers, Research Institute for Symbolic Computation, J. Kepler University Linz, Austria}
 \email{mkauers@risc.uni-linz.ac.at}
 \thanks{$^\ast$ supported in part by the Austrian Science Foundation (FWF)
  grants SFB F1305 and P19462-N18. This work was done during a visit of M.\,K.
  to D.\,Z. at Rutgers University.}

 \author[Doron Zeilberger]{Doron Zeilberger\,$^\dagger$}
 \address{Doron Zeilberger, Rutgers University, Piscataway, New Jersey, USA}
 \email{zeilberg@math.rutgers.edu}
 \thanks{$^\dagger$ supported in part by the USA National Science Foundation}

 \begin{abstract}
   We consider some multivariate rational functions
   which have (or are conjectured to have) only positive coefficients in their
   series expansion.
   We consider an operator that preserves positivity of series coefficients,
   and apply the inverse of this operator to the rational functions.
   We obtain new rational functions which seem to have only positive coefficients,
   whose positivity would imply positivity of the original series,
   and which, in a certain sense, cannot be improved any further.
 \end{abstract}

 \maketitle

 \section{Introduction}

 Are all the coefficients in the multivariate series expansion about the origin of
 \[
  \frac1{1-x-y-z-w+\tfrac23(xy+xz+xw+yz+yw+zw)}
 \]
 positive? Nobody knows. For a similar rational function in three variables, 
 Szeg\"o~\cite{szego33} has shown positivity of the series coefficients using 
 involved arguments. His dissatisfaction with the discrepancy between the
 simplicity of the statement and the sophistication of the methods he used
 in his proof has motivated further research about positivity of the
 series coefficients of multivariate rational functions.
 For several rational functions, including Szeg\"o's, 
 there are now simple proofs for the positivity of their coefficients available. 
 For others, including the one quoted above~\cite{askey72}, the positivity of their 
 coefficients are long-standing and still open conjectures. 

 In this paper, we consider the positivity problem in connection with the
 operator $T_p$ ($p\geq0$) defined as follows:
 \begin{alignat*}1
  &T_p\colon\set R[\![x_1,\dots,x_n]\!]\to\set R[\![x_1,\dots,x_n]\!],\\
  &(T_p f)(x_1,\dots,x_n):=
      \frac{f\Bigl(\frac{px_1}{1-(1-p)x_1},\dots,\frac{px_n}{1-(1-p)x_n}\Bigr)}
           {\scriptstyle(1-(1-p)x_1)\cdots(1-(1-p)x_n)}
 \end{alignat*}
 By construction, the operator $T_p$ preserves positive coefficients for any 
 $0\leq p\leq 1$, i.e., if a power series $f$ has positive coefficients, then the power
 series $T_pf$ has positive coefficients as well, for any $0\leq p\leq1$.  For example, via
 \[
  T_{1/2}\Bigl(\frac1{1-x-y-z+4xyz}\Bigr)=\frac1{1-x-y-z+\tfrac34(xy+xz+yz)},
 \]
 positivity of the former rational function~\cite{askey77} implies positivity of 
 Szeg\"o's rational function~\cite{szego33}. 
 This is a fortunate relation, because the positivity of the former can be shown
 directly by a simple argument~\cite{gillis79} while this is not as easily possible 
 for the latter~\cite{kauers07b}. (Straub~\cite{straub07} gives a different positivity 
 preserving operator also connecting these two functions.)

 This suggests applying the operator $T_p$ ``backwards'' to a rational function~$f$
 for which positivity of the coefficients is conjectured, in the hope that this leads
 to a rational function which again has positive coefficients, and for which positivity
 of the coefficients is easier to prove. 
 We present some empirical results in this direction. 
 Our results may or may not lead closer to rigorous proofs of some open problems.
 In either case, we also find them interesting in their own right.

 \section{Sharp Improvements}

 Given a rational function~$f$, we are interested in parameters $p\in[0,1]$ such that
 $T_p^{-1}f$ has positive series coefficients. Because of $T_p^{-1}=T_{1/p}$, this is
 equivalent to asking for parameters $p\geq1$ such that $T_pf$ has positive series
 coefficients.  Clearly, the set of all $p\geq0$ such that $T_pf$ has positive
 coefficients forms an interval $[0,p_{\max})$ with a characteristic upper bound
 $p_{\max}$ for each particular~$f$.
 Computer experiments have led to the following empirical results. 

 \begin{ethm}\label{thm:1}
   Let $f(x,y,z)=1/(1-x-y-z+4xyz)$.
   Let $p_0$ be the real root of $2x^3-3x^2-1$ with $p_0\approx1.68$.
   Then $p_0=p_{\max}$.
 \end{ethm}
 \begin{proof}[Evidence]
   \begin{enumerate}
   \item % checked.
     $p_{\max}$ cannot be larger than $p_0$, because the particular
     coefficient $\<xyz>T_pf=1+3p^2-2p^3$ fails to be positive for $p\geq p_0$.
   \item % checked.
     CAD computations confirm that all terms $\<x^ny^mz^k>T_pf$ with $0\leq n,m,k\leq50$
     are positive for any $0<p<p_0$.
   \item % checked. 
     For $p=2430275/1448618$, all terms $\<x^ny^mz^k>T_pf$ with $0\leq n,m,k\leq 100$ 
     are positive.
     This $p$ is the 15th convergent to $p_0$ and only about
     $10^{-14}$ smaller than this value. 
%    \item % OPEN.
%      \ [[FPSAC-like arguments for diagonals]]
%      not even the main diagonal admits a proof for p close to p_max
   \item % checked.
     For each specific choice of $m,k$, the terms $\<x^ny^mz^k>T_pf$ are polynomials
     in $n$ (and~$p$) of degree $m+k$ with respect to~$n$. 
     For $0\leq m,k\leq10$, CAD computations confirm that these are
     postive for all $n\geq1$ and all $0<p<p_0$.
   \end{enumerate}
 \end{proof}

 \begin{ethm}\label{thm:2}
   Let $f(x,y,z,w)=1/(1-x-y-z-w+\tfrac23(xy+xz+xw+yz+yw+zw))$.
   Let $p_0$ be the real root of $x^4-6x^2-3$ with $p_0\approx2.54$.
   Then $p_0=p_{\max}$.
 \end{ethm}
 \begin{proof}[Evidence]
   \begin{enumerate}
   \item % checked.
     $p_{\max}$ cannot be larger than $p_0$, because the particular
     coefficient $\<xyzw>T_pf=3+6p^2-p^4$ fails to be positive for $p\geq p_0$.
   \item % checked. 
     CAD computations confirm that all terms $\<x^ny^mz^kw^l>T_pf$ with $0\leq n,m,k,l\leq25$
     are positive for any $0<p<p_0$.
   \item % checked.
     For $p=730647/287378$, all terms $\<x^ny^mz^kw^l>T_pf$ with $0\leq n,m,k,l\leq240$
     are positive.
     This $p$ is the 15th convergent to $p_0$ and only about
     $10^{-12}$ smaller than this value. 
%    \item % OPEN.
%       \ [[FPSAC-like arguments for diagonals]]
%      not even the main diagonal admits a proof for p close to p_max
   \item % checked.
     For each specific choice of $m,k,l$, the terms $\<x^ny^mz^kw^l>T_pf$ are polynomials
     in $n$ (and~$p$) of degree $m+k+l$ with respect to~$n$. 
     For $0\leq m,k,l\leq 5$, CAD computations confirm that these polynomial are
     postive for all $n\geq1$ and all $0<p<p_0$.
   \end{enumerate}
 \end{proof}

 For the rational function~$f$ considered in Statement~\ref{thm:2},
 our hope was dispelled that a direct proof for the positivity of $T_{p_{\max}}f$ could 
 be found more easily than for~$f$ itself. However, some ``suboptimal'' values of $p$ do
 lead to rational functions which have a promising shape. For instance, we found that
 \begin{alignat*}1
   &T_{\sqrt3}\Bigl(\frac1{1-x-y-z-w+\tfrac23(xy+xz+xw+yz+yw+zw)}\Bigr)\\
   &\quad=\frac1{1-x-y-z-w+2(xyz+xyw+xzw+yzw)+4xyzw}.
 \end{alignat*}
 Also note that it seems to be a coincidence that $p_{\max}$ is determined by the 
 coefficient of $xyzw$ in $T_pf$, because this does not hold in the expansion of
 \[
   f(x,y,z,w)=\frac1{1-x-y-z-w+\tfrac{64}{27}(xyz+xyw+xzw+yzw)}
 \]
 which is conjectured to have positive coefficients~\cite{kauers07b}. 
 Here we have $p_{\max}<1.66$, by inspection of the coefficients $\<x^ny^mz^kw^l>T_pf$
 for $0\leq n,m,k,l\leq 100$, while $\<xyzw>T_pf=-\frac{13}{27}p^4-\frac{40}{27}p^3+6 p^2+1$ 
 is positive for $p<2.36$. 

 \section{Asymptotically Positive Coefficients}

 Inspection of initial coefficients of
 \[
  \frac1{1-x-y-z-w+\tfrac{64}{27}(xyz+xyw+xzw+yzw)}
 \]
 suggests values for $p_{\max}$ that become smaller and smaller as the amount of initial values
 taken into consideration increases. Is the ``real'' value $p_{\max}$ 
 determined by the asymptotic behaviour of the coefficients for general~$p$?

\begin{figure}
   \leavevmode\hfill
   a)~\epsfig{file=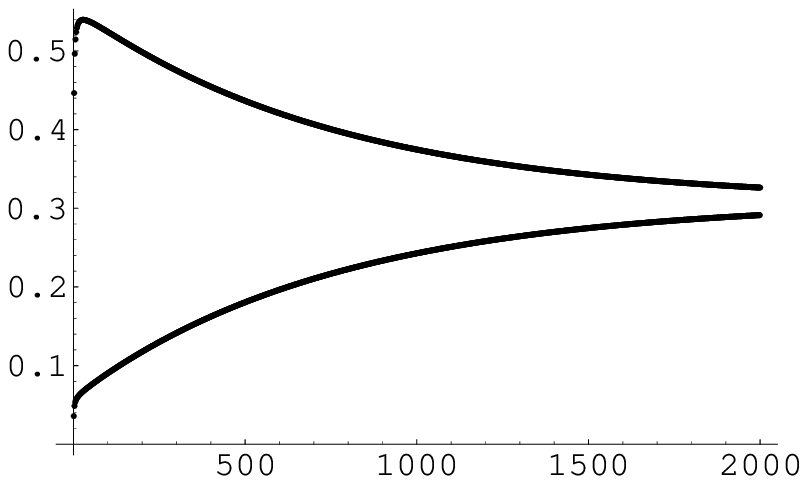,width=.33\hsize}\hfill
   b)~\epsfig{file=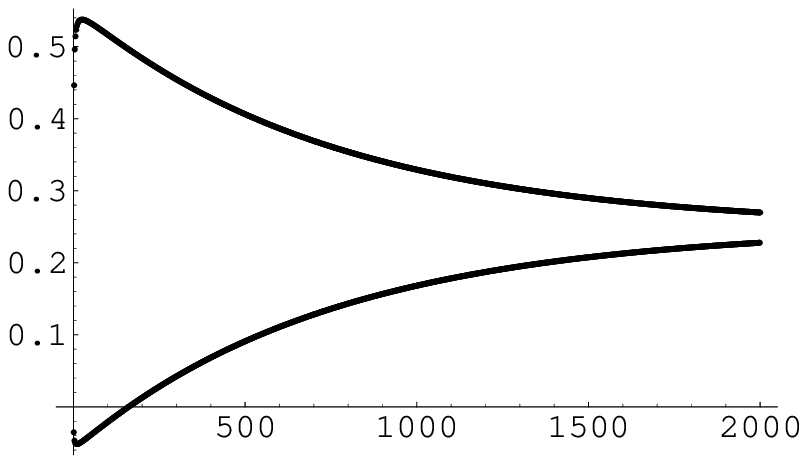,width=.33\hsize}\hfill\null
   \caption{}\label{fig:1}
 \end{figure}

 Clearly, it is hard to extract conjectures about asyptotic behaviour by just looking
 at initial values. Instead, such information is better extracted from suitable recurrence
 equations by looking at the characteristic polynomial and the indicial equation of
 the recurrence~\cite{wimp85}. 
 Using computer algebra, obtaining recurrence equations for the coefficient sequences is
 an easy task.
 Often, the asymptotics can be rigorously determined from a recurrence up 
 to a constant multiple~$K$, which cannot be determined exactly, but for which numeric
 approximations can be found. For instance, if $p>p_0:=(15+3 \sqrt{33})/2$, we have
 \begin{alignat*}1
   a_n&:=\<x^ny^nz^nw^0>T_p\Bigl(1/(1-x-y-z-w+\tfrac{64}{27}(xyz+xyw+xzw+yzw))\Bigr)\\
      &\qquad
      \sim K\Bigl(\frac{\left(155+27 \sqrt{33}\right) \left(-4 p+3
       \sqrt{33}-15\right)^3}{3456}\Bigr)^nn^{-1}
       \quad(n\to\infty)
 \end{alignat*}
 for $K\gtrsim0.291$ (Figure~\ref{fig:1}a shows $a_n\big/((\cdots)^nn^{-1})$ 
 for $p=p_0+\tfrac1{10}$, supporting the estimate for~$K$). This is oscillating. 
 For $1<p<p_0$, the asymptotic behaviour turns into
 \[
   a_n\sim K\bigl(1+\tfrac53p\bigr)^{3n}n^{-1}
   \quad(n\to\infty),
 \]
 for $K\gtrsim0.227$ (Figure~\ref{fig:1}b shows $a_n\big/((\cdots)^nn^{-1})$
 for $p=p_0-\tfrac1{10}$, supporting the estimate for~$K$). 
 This is not oscillating, but ultimately positive.
 This supports the conjecture $p_{\max}<p_0\approx16.1168$, 
 which is little news, however, as we already know $p_{\max}<1.66$ by inspection 
 of initial values.  Other paths to infinity that we tried 
 do not give sharper bounds on~$p_{\max}$.  
 So it seems that $p_{\max}$ in this example is determined
 neither by the initial coefficients, nor by the coefficients at infinity, 
 but somehow by the coefficients ``in the middle''.

 We can consider asymptotic positivity of coefficients as an independent question
 which may also be asked for the rational functions considered in Statements~\ref{thm:1}
 and~\ref{thm:2}: What are the values $p\geq p_{\max}\geq1$ such that the series
 coefficients of $T_pf$ are \emph{ultimately} positive?
 Denote by $p_{\max}^\infty$ the supremum of these parameters.
 We have carried out computer experiments in search for $p_{\max}^\infty$, and we
 obtained the following empirical results.

 \begin{ethm}
   Let $f(x,y,z)=1/(1-x-y-z+4xyz)$. Then $p_{\max}^\infty=2$.
 \end{ethm}
 
 \begin{proof}[Evidence]
   Let $\epsilon>0$ (sufficiently small) and $a_{n,m,k}:=\<x^ny^mz^k>T_{2-\epsilon}f$.
   \begin{enumerate}
   \item % checked. 
     First of all, we have $p_{\max}^\infty\leq2$, because for $\epsilon=0$, 
     the asymptotics on the main diagonal is
     \[
       a_{n,n,n}\sim K(-27)^nn^{-2/3}\quad(n\to\infty)
     \]
     with $K\gtrsim0.25$, i.e., $a_{n,n,n}$ is ultimately oscillating for $\epsilon=0$.
   \item % checked.
     Let $m,k\geq0$ be fixed and consider $a_{n,m,k}$ as a sequence in~$n$.
     A direct calculation shows that
     \begin{alignat*}1
       a_{n,m,k}&=\sum_{r=0}^m\sum_{t=0}^k\sum_{s=0}^t
         (-1)^{r+s}
         \binom nr\binom{n+m-r}{m-r}\binom{n+m-2r}s\binom{n+m-r+t-s}{t-s}\\
         &\qquad\times\binom{r}{k-t}
         (3-\epsilon)^{r+k-t+s}(3-2\epsilon)^{k-t}(\epsilon-1)^{r-k+t+s}\\
         &=\frac{(2-\epsilon)^{2(m+k)}}{m!k!}n^{m+k}+o(n^{m+k})\quad(n\to\infty),
     \end{alignat*}
     which is positive for $n\to\infty$.
   \item\label{item:4:3} 
     % checked.
     For arbitrary (symbolic) $i\geq0$ and the particular values $0\leq j\leq 3$, the sequence
     $a_{n,n+i,j}$ satisfies a recurrence equation of order~3 which gives rise to
     \begin{alignat*}3
       a_{n,n+i,j}&\sim K (3-\epsilon)^{2n} n^{-1/2}&&\quad(n\to\infty)
     \end{alignat*}
     for some constants~$K$ depending on $i,j$, and~$\epsilon$.
     % plots checked for 0<=i,j<=3, but for i=3,j=1 the plot-evidence is weak.
     Numeric computations suggest that these constants are positive, and hence,
     $a_{n,n+i,j}$ is positive for $n\to\infty$.
%
%      For instance, Fig.~\ref{fig:1} shows $a_{n,n+1,1}/\bigl((3-\epsilon)^{2n}/n^{1/2}\bigr)$ 
%      when $\epsilon=0.1$~(left) and when $\epsilon=0.01$~(right), 
%      indicating that the limit $K$ is positive in these cases. 
%      For different choices of $i,j,\epsilon$, the plots have a similar shape. 
   \item\label{item:4:4} % checked.
     For the particular values $0\leq i,j\leq 2$, 
     the sequence $a_{n,n+i,n+j}$ satisfies a recurrence equation of order~3
     which gives rise to 
     \[
       a_{n,n+i,n+j}\sim K(3-\epsilon)^{3n}n^{-1}\quad(n\to\infty)
     \]
     for some constants~$K$ depending on $i,j$, and~$\epsilon$. 
     % numerics checked for 0<=i,j<=2 and epsilon=.1
     Numeric computations suggest that these constants are positive, and hence,
     $a_{n,n+i,n+j}$ is positive for $n\to\infty$.
   \end{enumerate}
 \end{proof}

%  \begin{figure}
%    \leavevmode\hfill
%    \epsfig{file=fig1a.eps,width=.33\hsize}\hfill
%    \epsfig{file=fig1b.eps,width=.33\hsize}\hfill\null
%    \caption{Limiting behaviour of $\frac{a_{n,n+1,1}}{(3-\epsilon)^{2n}/\sqrt n}$ 
%       for $\epsilon=.1$ (left) and $\epsilon=.01$ (right).}
%    \label{fig:1}
%  \end{figure}

 In parts \ref{item:4:3} and~\ref{item:4:4}, we could not carry out the arguments
 for both $i$ and $j$ being generic. 
 We did find a recurrence equation of order~6 for $a_{n,n+i,n+j}$ for generic $i,j$,
 with polynomial coefficients of total degree~16 with respect to $n,i,j$, but this
 recurrence was way to big for further processing.  

 \begin{ethm}
   Let $f(x,y,z,w)=1/(1-x-y-z-w+\tfrac23(xy+xz+xw+yz+yw+zw))$. Then $p_{\max}^\infty=3$.
 \end{ethm}

\begin{figure}
   \leavevmode\null\hfill
   \epsfig{file=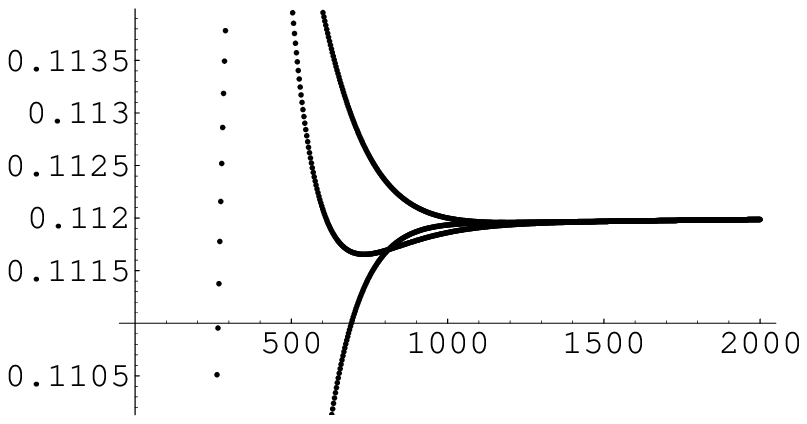,width=.33\hsize}\hfill\null
   \caption{}
   \label{fig:2}
 \end{figure}

 \begin{proof}[Evidence]
   Let $p\geq1$ and $a_{n,m,k,l}:=\<x^ny^mz^kw^l>T_pf$.
   \begin{enumerate}
   \item\label{item:4:1} % checked.
     First of all, we have $p_{\max}^\infty\leq3$ because for $p>3$, the asymptotics on
     the main diagonal is determined by the two complex conjugated roots
     % in fact, there is a real dominating root, but it is fake if p>3 and contributes
     % only for p<3, see below.
     \[
       \frac{9+30p^2-7p^4\pm 4p(p^2+3)\sqrt{6-2p^2}}{9}.
     \]
     Their modulus is $(p^2-1)^2$.
     As $(p^2-1)^2$ itself is not a characteristic root, it follows~\cite{gerhold06c}
     that $a_{n,n,n,n}$ is ultimately oscillating for $p>3$.     
   \item % OPEN.
     For $i,j,k\geq0$ fixed, $a_{n,i,j,k}$ is a polynomial in $n$ of degree $i+j+k+1$
     whose leading coefficient is $p^{2(i+j+k)}/3^{i+j+k}/i!/j!/k!$. Therefore
     $a_{n,i,j,k}$ is positive for $n\to\infty$ regardless of~$p$.
   \item % checked.
     For the particular values $i=0,1$ and $0\leq j,k\leq 2$, the sequence
     $a_{n,n+i,j,k}$ satisfies a recurrence equation of order 3 which gives rise to
     % rec has poly coeffs of deg i+j+k. charpoly of (i,j,k)=(1,0,1) has a fake root
     \[
       a_{n,n+i,j,k}\sim K\frac{(p+\sqrt3)^{2n}}{3^n}n^{j+k-\tfrac12}\quad(n\to\infty)
     \]
     for some constants $K$ depending on $i,j,k$, and~$p$.
     % checked.
     Numeric computations suggest that these constants are positive, and hence, 
     $a_{n,n+i,j,k}$ is positive for $n\to\infty$ regardless of~$p$.
   \item % checked.
     For the particular values $0\leq j,k\leq 1$, the sequence
     $a_{n,n,n+j,k}$ satisfies a recurrence equation of order~4 which gives rise to
     % the dominating charpoly root (9-p^4) must have coefficient zero, because
     % (9-p^4)^n does not match the degree 3n of a_{n,n,n,0} wrt. p. the next 
     % in order is (1+p)^3.
     \[
       a_{n,n,n+j,k}\sim K(1+p)^{3n}n^{-1}\quad(n\to\infty)
     \]
     for some constants~$K$ depending on $j,k$, and~$p$.
     % checked.
     Numeric computations suggest that these constants are positive, and hence,
     $a_{n,n,n+j,k}$ is positive for $n\to\infty$ regardless of~$p$.
   \item % checked.
     The main diagonal $a_{n,n,n,n}$ satisfies a recurrence of order 4 which gives
     rise to 
     \[
       a_{n,n,n,n}\sim K
       \Bigl(64+\tfrac{1}{27}(p^2-9)(2 p^4+9 p^2+189-2(p^2+3)^{3/2} p)\Bigr)^n
       n^{-3/2}
     \]
     for some constant~$K$ depending on~$p$.

     Numeric computations suggest that $K$ is positive for $p<3$.
     For example, Figure~\ref{fig:2} shows the quotients $a_{n,n,n,n}/((\cdots)^nn^{-3/2})$ for
     $p=3-\tfrac1{100}$.
   \end{enumerate}
 \end{proof}

%  The main diagonal $a_{n,n,n,n}$ satisfies a recurrence of order 4 with polynomial
%  coefficients of degree~7. This recurrence, however, does not admit asymptotic 
%  analysis owing to two complex conjugated dominant roots of the characteristic
%  polynomial (at least for $p>\sqrt3$). 

 Stronger evidence in support of the conjectures made in the paper is currently 
 beyond our computational and methodical capabilities. So are rigorous proofs.

 \bibliographystyle{plain}
 \bibliography{all}

\begin{thebibliography}{1}

\bibitem{askey72}
Richard Askey and George Gasper.
\newblock Certain rational functions whose power series have positive
  coefficients.
\newblock {\em The American Mathematical Monthly}, 79(4):327--341, 1972.

\bibitem{askey77}
Richard Askey and George Gasper.
\newblock Convolution structures for {L}aguerre polynomials.
\newblock {\em Journal d'Analyse Math.}, 31:46--68, 1977.

\bibitem{gerhold06c}
Stefan Gerhold.
\newblock Point lattices and oscillating recurrence sequences.
\newblock {\em Journal of Difference Equations and Applications},
  11(6):515--533, 2005.

\bibitem{gillis79}
Joseph Gillis and J.~Kleeman.
\newblock A combinatorial proof of a positivity result.
\newblock {\em Mathematical Proceedings of the Cambridge Philosopical Society},
  86:13--19, 1979.

\bibitem{kauers07b}
Manuel Kauers.
\newblock Computer algebra and power series with positive coefficients.
\newblock In {\em Proceedings of FPSAC'07}, 2007.

\bibitem{straub07}
Armin Straub.
\newblock Positivity of {S}zeg{\"o}s rational function.
\newblock {\em Advances in Applied Mathematics}, 2007.
\newblock to appear.

\bibitem{szego33}
Gabor Szeg{\"o}.
\newblock {\"U}ber gewisse {P}otenzreihen mit lauter positiven {K}oeffizienten.
\newblock {\em Mathematische Zeitschriften}, 37(1):674--688, 1933.

\bibitem{wimp85}
Jet Wimp and Doron Zeilberger.
\newblock Resurrecting the asymptotics of linear recurrences.
\newblock {\em Journal of Mathematical Analysis and Applications},
  111:162--176, 1985.

\end{thebibliography}
 
\end{document}